\documentclass[12pt]{article}
\usepackage{amssymb,amsmath,amsfonts,latexsym}
\begin{document}
\title{Multivariate truncated moments problems \\ and maximum entropy}
\author{C.-G. Ambrozie\footnote{Supported by the Czech grants IAA100190903 GAAV, 201/09/0473 GACR (RVO: 67985840) and the Romanian grant CNCS � UEFISCDI, no. PN-II-ID-PCE-2011-3-0119}}
\date{}

\newtheorem{remark}{Remark}
\newtheorem{remarks}[remark]{Remarks}
\newtheorem{proposition}[remark]{Proposition}
\newtheorem{examples}[remark]{Examples}
\newtheorem{example}[remark]{Example}
\newtheorem{lemma}[remark]{Lemma} 
\newtheorem{corollary}[remark]{Corollary}
\newtheorem{theorem}[remark]{Theorem}

\maketitle
\begin{abstract}
We characterize the existence of the Lebesgue integrable solutions of the truncated problem of moments in several  
variables on  unbounded  supports  by the existence of some maximum entropy -- type 
representing densities and discuss a few topics on their approximation in a particular case, of two variables and  4th order moments.
%Also we present a new method
 % of finding this density for  a certain class of such  problems, that we exemplify for $n=2$. 
%, and by the  boundedness of the corresponding Lagrangian functions.

Keywords: moments problem, representing measure, entropy

MSC-clas:  44A60 (Primary) 49J99 (Secondary)
\end{abstract}

\section{Introduction}

In this work we consider the problem of moments in the following context.
Let $T\subset \mathbb{R}^n$ be a closed subset,  where  $n\in \mathbb{N}$ is fixed.
Let $I\subset (\mathbb{Z}_{+})^n$ be finite such that $0\in I$, where $\mathbb{Z}_+ =\mathbb{N}\cup \{ 0\}$. 
Fix a set  $g=(g_i )_{i\in I}$  of real numbers $g_i $ with $g_0 =1$. 
The problem under consideration is to establish if there exist (classes of) 
Lebesgue measurable functions $f\geq 0$ a.e. (almost everywhere) on $T$, such that $\int_T |\, t^i |\, f(t)\, dt <\infty$ and
\begin{equation}
  \label{star}  
\int_T t^i f(t)\, dt =g_i \mbox{ }\mbox{ }\mbox{ }\mbox{ }\, (i\in I)
\end{equation}
and find such solutions $f$. 
As usual $dt=dt_1 \ldots dt_n$ and $t^i =t_{1}^{i_1} \cdots t_{n}^{i_n}$ for any multiindex $i=(i_1 ,\ldots ,i_n )\in \mathbb{Z}_{+}^n$ where
$t=(t_1 ,\ldots , t_n )$. 
In this case we call $f$  a {\it representing density} for $g$, and $g_i$ the {\it moments} of $f$. 
In general $T$ is unbounded and usually $I=\{ i:|i|\leq 2k\}$ for  $k\in \mathbb{N}$, where $|i|=i_1 +\cdots +i_n$. 

Generally  a problem of moments \cite{1},  \cite{19}, called also {\it $T$-problem of moments} when $T$ is given (R.E. Curto, L.A. Fialkow \cite{CuHr}),  is   concerned with the existence of
an arbitrary Borel measure $\mu \geq 0$ supported on $T$ such that $\int_T t^i d\mu (t)=g_i $ for $i\in I$, in which case one calls
$\mu$ a {\it representing measure} of $g$. 
The feasibility of (\ref{star})  characterizes  the dense interior of the convex cone of all data $g$ having representing measures,
provided that all $t\in T$ are density points  and $I$ is a union of intervals $[0,i]:=\{ j\in \mathbb{Z}_{+}^n :0\leq j_k \leq i_k ,\, 1\! \leq\! k\!\leq\! n\}$, see [Theorems 5, 6, \cite{RHav}] (and   M. Junk, [Theorem A.1, \cite{junk}]  in a slightly different context).
For our purpose here we only require that the Lebesgue measure of $T$
be $\not =0$.

Author's main contributions are contained in the statements 4 -- 9.  In particular, by  Corollary \ref{c2} for example, for each fixed $\epsilon >0$ we  characterize the feasibility of (\ref{star}) by the existence of a
(unique) $f_*$ minimizing  $\int_T f\ln f\, dt +\epsilon \int_T \| t\|^{2k+2}f(t)\, dt$ amongst all solutions, which is equivalent to the existence of
a (unique)
vector $\lambda^* =(\lambda_{i}^* )_{|i|\leq 2k}$ 
maximizing the associated  Lagrangian $L(\lambda )=L_\epsilon (\lambda )=\sum_{|i|\leq 2k}g_i \lambda_i -\int_T e^{\sum_{|i|\leq 2k}\lambda_i t^i -\epsilon \| t\|^{2k+2}} dt,$
 in which case $f_* (t)=e^{\sum_{|i|\leq 2k}\lambda_{i}^* t^i -\epsilon \| t\|^{2k+2}}$ where $\| t\| =(\sum_{j=1}^n t_{j}^2 )^{1/2}$.  The more general 
formulation of the main result Theorem \ref{lll}  aims to cover also other cases like $T= $compact   with $\epsilon =0$ \cite{Leww}.  

%Other results are Propositions  \ref{suplim}, \ref{p1} and Corollaries \ref{cc}, \ref{exempluu}. 

%In a particular case for $n=2$, we give a way of finding $\lambda^*$ (Proposition \ref{exempluu}).

Maximizing the Boltzmann-Shannon's {\it entropy}  $H(f)=-\int_{\mathcal{T}} f\ln f\, d\mu$ on a probability space $(\mathcal{T}, \mu )$ subject to various  restrictions $\int_{\mathcal{T}}a_i f d\mu \! =\! g_i$ ($i\in I$) is  a well-known principle  in   statistical mechanics 
and information 
theory 
 \cite{6},   \cite{12'}, \cite{junk}, \cite{MP}. The maximum of $H$
is attained on the  unbiased
probability distribution $f_*$ on a partial knowledge, of the prescribed  average values $g_i$ of some  random variables  \cite{Boltz}, \cite{6}, \cite{12'}. 
Typically $f_*$ is obtained  by maximizing a function $L$  (the Lagrangian) convex conjugate to $-H$  \cite{BL}, \cite{sharp}, \cite{Levermore}, \cite{Moreau}, 
 \cite{rod},   which leads to  characterizations ($\sup \! /\! \max  L\! <\! \infty $) of the feasibility of the primal problem -in our case (\ref{star}). 
One may  consider more general measures $\mu \!\geq \! 0$  or
functionals like $H( f) \! =\!  -\mbox{tr} (f\ln f)$, $\mbox{tr} (\ln f)$ where $f\! =$ positive definite matrix for  noncommutative moments 
 [Theorems 2,3, \cite{posmatr}], \cite{BWo}.

While the case  $T= $compact  was known long before (A.S. Lewis,  \cite{Leww}),
 the similar problems with unbounded support $T$ (or unbounded moments $a_i$) are usually  difficult and  still studied, see R.V. Abramov \cite{abr}, 
J.M. Borwein \cite{Boltz},  M. Junk \cite{junk}, C.D. Hauck, C.D. Levermore, A.L. Tits \cite{Levermore} and others   \cite{junkl}, \cite{leo}, \cite{gmc}, \cite{junke}, \cite{tag}. We mention that the feasibility of  (\ref{star}) has been  characterized  in G. Blekherman, J.B. Lasserre \cite{bla}, avoiding entropy maximization but also in  Lagrangian terms.
If  $T$ is unbounded,  Corollary \ref{c2} cannot be improved to $\epsilon    \! =\!  0$: there are examples of realizable, but degenerate data $g$  such that  the constrained $H$-maximization fails for $(\mathcal{T}, \mu )\! =\! (\mathbb{R}^n ,dt)$, see M. Junk and co-workers \cite{junkl}, \cite{Levermore}. For $H(f)\! =\! -\int_T f \ln f\, dt$, the  maximization of $L(\lambda )$ ($\! =\! L_0 (\lambda )\in [-\infty ,\infty )$) always holds, at a unique point $\lambda^*$ - using  for instance [Corollary 2.6, \cite{BL}], see also  \cite{junk}, \cite{Levermore}.  It follows,  by means of Fatou's lemma, for $I\! =\! \{ i\! :\!  |i|\! \leq\!  2k\}$, that $|t^i |e^{\sum_{|j|\leq 2k} \lambda^{*}_j t^j}\! \in\!  L^1 (T, dt)$  for all $|i|\! \leq\!  2k$,  $ \int_T t^i e^{\sum_{|j|\leq 2k} \lambda^{*}_j t^j}dt\! =\! g_i$ ($|i|\! <\! 2k$) but the equality may fail for $|i|\! =\! 2k$.
Namely the dual attainment $\sup L\! =\! \max L$  does hold, but  primal attainment $\sup_{f\in (\ref{star})}H(f)\! =\! \max_{f\in (\ref{star})}H(f)$ is also a  difficult topic  if  $a_i (t)$ (for instance $t^i $) are not in the dual of $L^1 (T)$. For these matters we refer the reader to   \cite{junk},  \cite{junkl}, \cite{Levermore}, \cite{MP} and for  applications to Boltzmann equations we mention also \cite{cerc}, \cite{junke},  \cite{au}, \cite{tag},  \cite{gmc}.

Originated in works by Stieltjes, Hausdorff,  Hamburger and Riesz,  the area of moments problems  saw extensive development in many directions, that we do not attempt to cover.  
There exist also other  approaches   to the multivariate   moments problems,
by  operator theoretic or convexity methods \cite{8}, \cite{11},  \cite{h},  \cite{14''}, \cite{SS}, \cite{20},  in particular
 a  truncated version of  Riesz-Haviland's theorem \cite{CuHr},  see also   \cite{KM},  \cite{PS}   for other results,
  related to   sums-of-squares representations of positive polynomials or
  polynomial optimization theory. These interesting topics are beyond the  goal of the present paper, that is focused on the  $H\, /\, L\, $-- maximization. 
\vspace{3 mm}

 {\bf Acknowledgements} The present work was supported by the Czech grants IAA 100190903 GA AV and 201/09/0473 GA CR (RVO: 67985840) and the Romanian grant CNCS � UEFISCDI, no. PN-II-ID-PCE-2011-3-0119.

I wish to express my  thanks to  professor Marian Fabian for drawing the results of the Fenchel duality theory to my attention. Also, I  thank professor Mihai Putinar for several interesting suggestions and relevant references. 

The author is also  indebted to the reviewers for their expert indications  and reference  to more recent works, that definitely improved the paper.
%\vspace{13 mm}

\section{Main results}
 
Fix $T$, $I$ and $g$ as stated in the Introduction. 
For any  measurable space $\mathcal{T}$ endowed with a $\sigma$-finite measure $\mu \geq 0$ and $1\leq p\leq \infty$,
the notations $L^p (\mathcal{T},\mu )$, $L_{+}^p (\mathcal{T},\mu )$ (sometimes,  $L^p (\mu )$, $L_{+}^p (\mu )$) have the usual meaning. 
We repeat below an  argument from  [Theorem 2.9, \cite{BL}],  adapted to our case. 
\begin{lemma}
\label{aproximare}
{\em (see \cite{BL})} Let $\mu \geq 0$ be a finite measure on $\mathcal{T}$.
 %be a measurable space endowed with a finite measure $\mu \geq 0$. 
 Let $x \in L^{1}_+ ( \mu )\setminus \{ 0\}$, and
$a_i \in L^1 (\mu )$ ($i\in I$) be a finite set of functions  such that  $\int_\mathcal{T} |a_i |x d\mu <\infty$ for all $i$ and $(a_i )_i$ are linearly independent on any subset of positive measure.
Then there is a sequence $(y_k )_{k\geq k_0} \subset L^\infty (\mu )$
such that  $x_k :=\min (x,k)+y_k \geq 0$ a.e., $\int_\mathcal{T} a_i x_k d\mu =\int_\mathcal{T} a_i x\, d\mu $ for all $i\in I$, $|y_k |\leq x$ and $y_k \to 0$ a.e. % as $k\to \infty$.

%------------------------------------------

\end{lemma}

{\it Proof}.  
Set $z_k \! =\! \min \! (x,k)$ for $k\! \geq\! 1$. Using $\{  \! x\! >\! 0\! \} \! =\! \cup_{l\geq 1}\{ \! x\! \geq \! 1/l\! \}$, we find a $\delta \in (0,1)$ and $\mathcal{T}_* \subset \mathcal{T}$ with  $\mu (\mathcal{T}_* )>0$ such that  $x(t)\geq \delta$ a.e.  on $\mathcal{T}_*$. %Endow $\mathbb{R}^N$ ($N=\mbox{card}\, I$) with a norm.
The linear  map $A:L^\infty (\mathcal{T}_* )\to \mathbb{R}^N$ ($N=\mbox{card}\, I$),  $Ay=( \int_{\mathcal{T}_*} a_i y\, d\mu )_{i} $  is surjective for  otherwise there is a  $(\lambda_i )_i \not =0$ orthogonal to its range, such that  $\sum_i \lambda_i \int_{\mathcal{T}_*}a_i y \, d\mu =0$ $\forall \, y$,
whence $\sum_i \lambda_i a_i =0$ a.e. on $\mathcal{T}_*$ that  is impossible. % because
%  the set of zeroes of a polynomial$\, \not =0$ has measure zero.
   Since $A$ has closed range,  there is a $c$ such that   $\inf_{w\in \ker\, A} \| y -w\|_\infty \leq c \| Ay \| $ $\forall$ $y\in L^\infty (\mathcal{T}_* )$.
   % is finite and $0$ is  interior  to $C=\{ Ay: y\in L^\infty (T_* ),\, \| y\|_\infty <\delta /2\}$.
   % and the reduced minimum modulus of $A$ is $>0$. 
By Lebesgue's theorem of dominated convergence, $\lim_k \int_\mathcal{T} a_i z_k d\mu =\int_\mathcal{T} a_i x\, d\mu $ for all $i$.
%Hence for  large $k$, $(\int_T a_i (x-z_k )\, d\mu )_i \in C$. %Moreover
 There are $y_k \in L^\infty (\mathcal{T})$  with $\mbox{supp}\, y_k \subset \mathcal{T}_*$
 such that  $\int_\mathcal{T} a_i y_k d\mu =\int_\mathcal{T} a_i (x-z_k )\, d\mu$ and, since  $ Ay_k  \to 0$,  we can choose them such that 
$\| y_k \|_\infty \to 0$. For large $k$,  $\| y_k \|_\infty \leq \delta /2$. On $\mathcal{T}_*$, $x\geq \min (x,k)\geq \delta >\delta /2\geq | y_k |$. Then $x_k \geq 0$ a.e.
%,
%while outside $T_*$, $y_k =0$. 
$\Box$
\vspace{3 mm}

Fenchel duality  
 deals with  minimizing  convex functions $\varphi \! :\! X\! \to\!  (-\infty ,\infty \, ]$ over convex subsets of locally convex spaces $X$, in connection with the dual problem of maximizing  $-\varphi^*$ where $\varphi^* $  is the {\it convex conjugate} of $\varphi$, called also its {\it Legendre-Fenchel transform}  \cite{BL}, \cite{rod}, \cite{Moreau}, \cite{sharp},  \cite{Levermore}; $\varphi$ must be 
{\it proper} ($\varphi\not\equiv \infty $). Letting the {\it effective domain} of $\varphi$ be ${\rm dom}\, \varphi \! =\! \{ x\in X: \varphi(x)\! <\! \infty\}$, $\varphi^*$ is defined on the dual of $X$ by 
 $\varphi^* (x^* )\! =\! \sup\{ \langle x,x^* \rangle \! -\! \varphi (x):x\in {\rm dom}\, \varphi  \} $.
 Typically, $\inf \varphi \! =\! \sup (-\varphi^*)$. 
Briefly speaking,
we set $\varphi(x)=-H (x)$ if $x\geq 0$ satisfies the equations of moments,  and $\varphi (x)=+\infty$ outside the set of  solutions. Then $\varphi$ is convex conjugate to $\varphi^* (  x^{* }) = \ln \int_T e^{\sum_{i}\lambda_i a_i } \, d\mu -\sum_{i}g_i \lambda_i $ for  $x^{*} =\sum_{i}\lambda_i a_i $, and $\varphi^* (x^* )=+\infty$ otherwise. Thus ${\rm dom}\, \varphi^*$ is the linear span of the $a_i$'s and (if $a_0 \equiv 1$) the   Lagrangian  $l:=-\varphi^* |_{{\rm dom}\, \varphi^*}$ is given by
 $\lambda \mapsto \!  -\ln \int_T e^{\sum_{i\in I\setminus \{ 0\}}\lambda_i (a_i -g_i )}d\mu $. Maximizing $l $ or $L$
 %=\sum_{i\in I}g_i \lambda_i -\int_T e^{\sum_{i\in I}\lambda_i a_i -1}d\mu$ 
  are equivalent problems. 
  We  rely on  J.M. Borwein and A.S. Lewis' results    \cite{BL}  concerned with $L$, providing dual attainment in a point $\lambda^*$. The equality $\inf \varphi=\sup (-\varphi^* )$  becomes here $P=D$.
 % (in a point $x^{*}_\lambda$ for a $\lambda =\lambda^*$). % $\inf F=\max (-F^* )=(-F^* (x^* ))$. 
 Although under different hypotheses, $L$ is analogous to the dual function $\psi$ from C.D. Hauck, C.D. Levermore and A.L. Tits [Section 4.1, \cite{Levermore}],  and
 would fit the case when ${\rm dom}\, L \cap \partial (  {\rm dom}\, L)\! =\! \emptyset$ in  M. Junk \cite{junk} 
except we do not have here a distinguished moment $a_m$ such that  $\lim\limits_{\| t\| \to \infty}\frac{|a_i (t)|}{1+a_m (t)} =0\,  (i\not \! =\! m)$.
% is  an infinite dimensional  Fenchel theoretic result similar to Theorem \ref{roo}.  %reminiscent to the Lagrange multipliers method on convex cones. 

The following J.M. Borwein and A.S. Lewis' result  from \cite{BL}  is the main Fenchel theoretic tool to be used later on in the proof of our Theorem \ref{lll}. 
% and to  paper \cite{BL} for an infinite dimensional result we need. T
\begin{theorem}
\label{maxim}
{\em [Corollary 2.6,\cite{BL}]}
Let $\mathcal{T}$ be a space with finite measure $\mu \geq 0$, $1\leq p\leq\infty$ and $a_i \in L^{q} (\mu ) $, $g_i \in \mathbb{R}$ for $i\in I$ ($=\, $finite) where $\frac{1}{p}+\frac{1}{q}=1$. Let  $\phi :\mathbb{R}\to (-\infty ,\infty ]$ be proper, convex and lower semicontinuous, with $(0,\infty )\subset \mbox{\em dom}\, \phi $.
Suppose there exist $x\in L^p (\mu  )$ with $x(t)> 0$ a.e. such that $\phi \circ x\in L^1 (\mu  )$ and
$ \int_\mathcal{T} a_i (t)\, x(t) \, d\mu (t)=g_i $ for $i\in I$. Then the values $P\in [-\infty , \infty )$ and $D\in [-\infty , \infty \, ]$ defined respectively by
$$
P\! =\!  \inf \{ \! \int_\mathcal{T} \!  \phi (x(t))\, d\mu (t)  :  x\in\! L^p (\mu ),\,   x\geq 0  \mbox{ a.e.},\mbox{ }
\phi \!\circ\! x\in\! L^1 (\mu  ), \! \int_\mathcal{T} \! \!  a_i xd\mu \! =\! g_i \,  \forall i \, \} 
$$
and
$$
D\! =\! \max  \{ \sum_{i\in I} g_i \lambda_i -  \int_\mathcal{T}  \! \phi^*   (   \sum_{i\in I} \lambda_i a_i (t)) \, d\mu (t) :
  \lambda_i  \in \mathbb{R}, \, \phi^* \! \! \circ\!  \sum_{i\in I} \lambda_i a_i \in L^1 (\mu ) \, \}
$$
are equal, $-\infty \leq P=D<\infty $ and the maximum  $D$ is attained.  
\end{theorem}

 \begin{remarks}
\label{conj} 

{\em  (a) Let $\phi$ be defined by $\phi (x)=x\ln x $ for $x>0$, $\phi (0)=0$ and $\phi (x)=+\infty$ for $x<0$. Then $\phi$ is proper, convex, lower semicontinuous, bounded from below, with effective domain $[0,\infty )$ and its convex conjugate is $\phi^* (y)=e^{y-1}$ for all $y\in
 \mathbb{R}$;  use to this aim that $\phi^* (y)=\sup_{x\geq 0}  (xy -x\ln x) $.}

 {\em (b) For the integrand $\phi $ defined at (a)  and $(\lambda_i )_{i\in I} =0$,  the constant function $(\phi^* \! \circ\! \sum_{i\in I} \lambda_i a_i )(t)\equiv\phi^* (0)$ is in $  L^1 (\mu )$. Thus  for  any data $a_i$, $g_i$ verifying the hypotheses of
Theorem \ref{maxim}, we obtain that          $-\infty <P=D<\infty$.}

{\em (c) Let $x \in L_{+}^1 (\mu )$ with $x\ln x\in L^1 (\mu)$ and $y_k \in L^1 (\mu )$ ($k\geq 1$) such that  $x_k :=\min (x,k) +y_k \geq 0$ a.e., 
$|y_k |\leq  x$   and $ y_k  \to 0$ a.e. as $k\to \infty$. % such that $x_k \geq 0$ a.e. 
By Lebesgue's dominated convergence theorem,  $\lim_{k} \int x_k \ln x_k d\mu =\int x\ln x\, d\mu$, since on $\{ t\! :\! x_k (t)\! \geq\! 1\}$,
$x_k \!\leq\! 2x \Rightarrow |x_k \ln x_k |\!\leq\! |2x \ln (2x)|$ while on $\{ t\! :\! x_k (t)\! <\! 1\}$, $|x_k \ln x_k |\!\leq\! 1/e$; hence $|x_k \ln x_k |\!\leq\! |2 x\ln x +(2\ln 2 )x| +1/e\in L^1 (\mu )$. 
}
\end{remarks}

%For example $a(t)\equiv 1$, $a(t)=(\| t\|^2 +1)^k$ and $a(t)=\| t\|^{2k}+1$ ($k\in \mathbb{N}$) satisfy this condition.

In Theorem \ref{lll}  the choice of the norm  on $\mathbb{R}^n$ is unimportant.
 We call a  function $a$  on $T$ {\em  independent of} $(t^i )_{i\in I\setminus \{ 0\}}$ if there are no  subsets $Z \!\subset \! T$ of positive  measure and constants $(c_i )_{i\in I\setminus \{ 0\}}$ such that $a\! =\! \sum_{i\in I\setminus\{ 0\}}c_i t^i$ on $ Z$.

\begin{theorem}
\label{lll}Let $T\! \subset\! \mathbb{R}^n$ be  closed, $I\!\subset\! \mathbb{Z}_{+}^n$  finite, $0\!\in\! I$ and $g\! =\! (g_i )_{i\in I}$  a set of numbers with $g_0 \! =\! 1$. Set  $m\! =\! \mbox{\em max}_{i\in I} |i|$. 
Let $a, \rho $ be measurable functions on $T$, $0\! <\! a, \rho \! <\! \infty $ a.e. such that   $\int_T e^{\frac{\| t\|^m +1}{\alpha \, a(t)}}\rho (t) \, dt\! <\! \infty$ for all $\alpha \! >\! 0$,
and $a$ is  independent of $(t^i )_{i\in I\setminus\{ 0\}}$.
The  statements {\em (a),  (b), (c)} are equivalent:

{\em (a)} There exist functions $f\in L^{1}_{+} (T,dt)$ such that $\int_T |t^i |f(t)dt<\infty$ and
\begin{equation}\label{problema}
\int_T t^i f(t)\, dt =g_i \mbox{ }\mbox{ }\mbox{ }\mbox{ }\, (i\in I);
\end{equation}

{\em (b)} There exists a particular solution $f_* \in L_{+}^1 (T, dt)$ of problem {\em (\ref{problema})}, maximizing the entropy functional $H\! =\! H_{\rho ,a}:L_{+}^1 (T,dt)\to [-\infty ,\infty )$ given by $$H(f)=-\int_T (\, \frac{af}{\rho}\ln \frac{af}{\rho}\, )\mbox{ } \rho \, dt$$
amongst all solutions;

{\em (c)} The  Lagrangian function $L=L_{\rho ,a,g }:\mathbb{R}^N \to [-\infty ,\infty )$ ($N=\mbox{\em card}\, I$) associated to the functional $H$ and the equations {\em (\ref{problema})}, given by $$L(\lambda )=\sum_{i\in I}g_i \lambda_i -\int_T e^{\sum_{i\in I}\lambda_i t^i /a(t)\mbox{ } -1}\, \rho (t) \, dt\mbox{ }\mbox{ }\mbox{ }\mbox{ }\mbox{ }\,
(\lambda =(\lambda_i )_{i\in I}),$$
 is bounded from above and attains its supremum in a point $\lambda^* =(\lambda_{i}^* )_{i\in I}$.
% of $\, \mathbb{R}^N$.

In this case $f_*$ and $\lambda^*$ are uniquely determined, $-H(f_* )=L(\lambda^* )$ and $$f_* (t)=a(t)^{-1}\rho (t)\, e^{\sum_{i\in I}\lambda_{i}^* t^i /a(t)-1}   \mbox{ }\mbox{ }\mbox{ }\mbox{ }\mbox{ }\mbox{ }
(t\in T),$$ 
in particular $H\not \equiv -\infty$ on the set of all solutions of {\em (\ref{problema})}, and
$$ 
\int_T \, t^i \, e^{\frac{1}{a(t)}\sum_{j\in I}\lambda_{j}^* t^j -1} a(t)^{-1}\rho (t) \, dt =g_i \mbox{ }\mbox{ }\mbox{ }\mbox{ }\, (i\in I).
$$ \end{theorem}

{\it Proof}.   Let 
   $a_i (t)=t^i /a(t)$ for $i\in I$ and $t\in T$. The condition on  $\rho$ and $a$  shows  that
  the measure $\mu :=\rho  \, dt$ on $T$ is finite  and, by means of 
the inequalities: $|t_j |\leq \| t\|$ ($:=(\sum_{j=1}^n t_{j}^2 )^{1/2}$) for $1\! \leq \! j\! \leq \! n$,
\begin{equation}
\label{mon}
|t^i |=|t_{1}^{i_1} \cdots t_{n}^{i_n }|\leq \| t\|_{}^{i_1 +\cdots +i_n}\leq (\| t\|_{}^{m }+1)^{|i |/m}\leq \| t\|_{}^{m}+1
\end{equation}
and $\sum_{i\in I}\lambda_i a_i (t)\leq \sum_{i\in I}|\lambda_i |\cdot \frac{\| t\|^m +1}{a(t)}$,
 that for every $\lambda =(\lambda_i )_{i\in I}\in \mathbb{R}^N $ 
\begin{equation}\label{ef}
g(\lambda ):= \int_T e^{\sum_{i\in I}\lambda_i a_i (t) -1}d\mu (t) <\infty .
\end{equation} 
By writing   $se^{s/\alpha} \leq e^{\beta s/\alpha}$ for large $\beta $ ($\geq \alpha /e +1$) and 
$s:=(\| t\|^m +1)/a(t)$,  
\begin{equation}\label{shift}
\int_T \frac{\| t\|^m +1}{a(t)}e^{\frac{\| t\|^m +1}{\alpha a(t)}}\, d\mu (t)<\infty\mbox{ }\mbox{ }\mbox{ }\, (\alpha >0).\end{equation}  
Then for every $\lambda =(\lambda_i )_{i\in I}$, by the inequalities (\ref{mon}) again, %we obtain as well
\begin{equation}\label{preven}
\int_T (\| t\|^m +1)a(t)^{-1}e^{\sum_{i\in I}\lambda_i a_i (t) -1}d\mu (t)<\infty .
\end{equation}
%(using $|i|\leq m$ for $i\in I$) 
%By (\ref{preven}) and (\ref{mon})
%we also  obtain  
Hence $\int_T a_i (t)e^{\sum_{i\in I}\lambda_i a_i (t)-1}d\mu (t)<\infty$, in particular $a_i \in L^1 (T,\mu )$ for $i\in I$.
Any of the statements (a) -- (c) implies that the Lebesgue measure of $T$ is strictly positive (finite or not), due to the condition $g_0 =1$.
 Then for every $f\in L_{+}^1 (T,dt)$, by Jensen's inequality for the 
 function $\phi (x):=x\ln x$ ($x\geq 0$),
 % and the probability measure $\frac{1}{\mu (T)}\mu$,
 $$
 H(f)=-\mu (T)\int_T \phi \, (\, \frac{af}{\rho }\, )\, \frac{d\mu }{\mu (T)}  \leq -\mu (T)\, \phi\,  (\int_T \frac{af}{\rho} \frac{d\mu }{\mu (T)} \, )\leq \mu (T)/e<\infty . 
 $$ 

 (a) $\Rightarrow$ (c).  Suppose that problem (\ref{problema}) has a solution $f$.
 The function $x :=af /\rho$ then  satisfies $\int_T | a_i |x d\mu <\infty $ and $\int_T a_i x d\mu =g_i$ for $i\in I$.
 By the original version [Theorem 2.9, \cite{BL}] of Lemma \ref{aproximare} (if   $x_k \! =\! \max (x,k)\!  +\! \frac{1}{k}\! +\! y_k$), there are  functions $\tilde{x} \in L_{}^\infty (T,\mu )$, $\tilde{x} > 0$ $\mu$-a.e. on $T$, such that $\int_T a_i (t) \tilde{x}(t)d\mu =g_i $ ($i\in I$). %in the original version \cite{BL} of Lemma \ref{aproximare}.  
 Here  $L^\infty (T,\mu )=L^\infty (T,dt)$ since  $\mu $ is equivalent to $dt$ on $T$.
For such $\tilde{x}$, the function $\phi \circ \tilde{x} =\tilde{x}\ln \tilde{x}$  belongs to $L^{\infty} (T )$, and hence, to $  L^{1}(T, \mu )$.
Then  we can use  Theorem \ref{maxim} for $\phi (x)=x\ln x$   and $p=\infty$,
see   Remark \ref{conj}, (a).  Let $P=\inf_x \int_T x\ln x\,  d\mu $ over the set of all  $x\in L^{\infty }_+(T)$  such that 
  \begin{equation}\label{adev}\int_T a_i x\, d\mu =g_i ,\mbox{ }\,
   i\in I\end{equation}and $D:=\sup L$. Then $-\infty <P=D<\infty$ with attainment in the dual problem, see  Remark \ref{conj}, (b). Therefore, 
$L(\lambda )=\sum_{i\in I}g_i \lambda_i -\int_T e^{\sum_{i\in I}\lambda_i a_i (t)-1}d\mu (t)$ is bounded from above on $\mathbb{R}^N$ and its supremum $D$  is attained. % in some point $\lambda^* =(\lambda_{i}^* )_{i\in I}$.
%there exists 

(c) $\Rightarrow$ (b). Assume there is a  $\lambda^* $ such that  $L(\lambda^* )=\max_{} L$.  As expected, we will derivate under the integral to show that
$x_* (t) :=e^{\sum_{i\in I}\lambda_{i}^* a_i (t) -1}$  satisfies (\ref{adev}) and moreover  maximizes $H_\mu (x):=-\int_T x\ln x\,  d\mu $ amongst all solutions from $L_{+}^1 (T,\mu )$.
Firstly, by (\ref{ef}), $\int_T x_* (t)d\mu (t)=g(\lambda^* ).$  By (\ref{mon}) and (\ref{preven}), $\int_T |a_i |x_* d\mu <\infty$ ($i\in I$).  
For any $\lambda $ we have    $L(\lambda )\leq L(\lambda^* )$, that is, by  (\ref{ef}),
%\sum_{i\in I}\tilde{\lambda}_i \tilde{\lambda}_i '$.
%Then for every $\tilde{\lambda}$  we have the inequality
\begin{equation}
\label{fff}
g(\lambda^* )\leq g(\lambda )+\sum_{i\in I}g_i  (\lambda_{i}^* -\lambda_i ).
\end{equation}
    Fix $j \in I$,  let $\varphi (t)=\pm\,  a_{j}(t)$
  and set $v=(v_i )_{i\in I}$ where $v_i =\pm \delta_{i\,  j}=\, $Kronecker's symbol (the signs agree). For any $\varepsilon >0$, set $\lambda_\varepsilon =\lambda^*
   +\varepsilon v$, namely $\lambda_\varepsilon =(\lambda_{\varepsilon i})_{i\in I}$ where $\lambda_{\varepsilon j}=\lambda_{j}^* \pm \varepsilon$ and $\lambda_{\varepsilon i}=\lambda_{i}^*$ for $i\not =j$.
Let  $F_\varepsilon (t)=\frac{1}{\varepsilon} x_* (t)(1-e^{\varepsilon \varphi (t)}).$  Note that
%\begin{equation}\label{ww}
\begin{equation}\label{limita}
\lim_{\varepsilon \to 0}F_\varepsilon (t)=-\varphi (t)x_* (t)\end{equation} and
$x_* e^{\varepsilon \varphi}=e^{\sum_{i\in I}\lambda_{i}^* a_i -1 +\varepsilon (\pm a_{j}) }=e^{\sum_{i\in I}\lambda_{\varepsilon i}a_i -1}$.
Then by (\ref{ef}) and (\ref{fff}),
%\end{equation} and \begin{equation}\label{www}
\begin{equation}\label{incauna}\int_T F_\varepsilon (t)d\mu (t)
= \frac{g(\lambda^* )-g(\lambda_\varepsilon )}{\varepsilon}\leq \mp g_{j}.\end{equation}
%\end{equation} 
By the estimates (\ref{mon}),
we may let  $y=\varphi (t)$ and $z=(\| t\|^{m} +1)/a(t)$  in the   inequality:
$
e^{- z}\frac{1-e^{\varepsilon y}}{\varepsilon}\geq -|y|
$ where $ z>0$, $y$ is real,   $|y|\leq z$ and $\varepsilon <1$.
Hence $
F_\varepsilon (t)\geq -x_* (t) |\varphi (t)|\cdot e^{(\| t\|^{m} +1)/a(t)}
$.  
 %the function $x_* (t) |\varphi (t)|\cdot e^{(\| t\|^{m} +1)/a(t)}$ in 
The right hand side is in $L^1 (T,\mu )$ 
by
the estimates: $|\varphi (t)|\leq (\| t\|^m +1)/a(t)$,
 $x_* (t)\leq e^{c(\| t\|^m +1)/a(t)}$  for a constant $c=c(\lambda^* )$, and  (\ref{shift}). 
Then we may apply Fatou's lemma for a sequence $\varepsilon =\varepsilon_k \to 0$ to obtain, by (\ref{limita}) and (\ref{incauna}), that
$$
\mp \int_T  a_{j}  x_*  d\mu  =-\int_T \varphi  x_*  d\mu  =\int_T
\lim_{\varepsilon \to 0}F_\varepsilon  d\mu  \leq 
 \liminf_{\varepsilon \to 0}\int_T F_\varepsilon  d\mu  \leq \mp g_{j}.
$$
%see also (\ref{fff}), (\ref{ww}) and (\ref{www}). Comparing the inequalities with $+$ and $-$ from above, 
Hence $\int_T a_{j} x_*  d\mu  =g_{j}$.  Since $j$ was arbitrary in $I$, $x_*$ is a solution of  (\ref{adev}). The  function $f_* :=\rho x_* /a$ is then a solution of (\ref{problema}). By 
(\ref{ef}) and (\ref{preven}), $x_* \ln x_* \in L^1 (T, \mu )$, i.e. $(af_* /\rho )\ln (af_* /\rho ) \in L^1 (T,\rho dt )$. Hence there are solutions $f$ of
 (\ref{problema}) such that  $H(f)> -\infty$.
By  the correspondence $f\leftrightarrow x=af/\rho$,  the fact that $f_*$ maximizes the functional $H$ given at (b)
is equivalent to saying
that
$\int_T x_* \ln x_* d\mu \leq \int_T x\ln x\, d\mu$ for all the solutions $x\in L^{1}_{+} (T,\mu ) $ of the 
problem (\ref{adev}). By Lemma \ref{aproximare} and Remark \ref{conj}, (c)
it suffices to show that $\int_T x_* \ln x_* d\mu \leq \int_T x\ln x\, d\mu  $
for any solution $x\in L^{\infty}_{+} (T)$ of (\ref{adev}). This holds  by
%. This follows by the usual estimates
%solutions $x$ from $ L^\infty (T)$, the existence of which was established in \cite.
$$
\int_T x\ln x\, d\mu  \geq P =D =\sum_{}\lambda_{i}^* g_i -\int_T e^{{\sum_{i}\lambda_{i}^* a_i -1}}d\mu =\int_T x_* \ln x_* d\mu
.$$

%---------------------------------------------------
%The remaining implications are  trivial. 

The conclusion $P=D$ of Theorem \ref{maxim} provides  $-H(f_* )=L(\lambda^* )$.
The uniqueness of $\lambda^*$ and $f_*$ (or, equivalently, $x_*$) follow from the strict convexity of $-L$, resp. $-H_\mu $ and the fact that $T$ is not negligible,  whence $p|_T =0 \mbox{ a.e.}  \Rightarrow p=0$ for any polynomial $p=\sum_{i\in I}\lambda_{i} X^i$  (the
zeroes sets   of  nonconstant polynomials are algebraic varieties, and so have null Lebesgue measure). $\Box$
\vspace{3 mm}

%In principle, for small $n$ and suitable choices of $\rho ,a$ we can obtain a vector $\lambda^*$, and so a solution $f^*$  by numerically
%maximizing  $L$. % seems to be possible for small $n$.

%In rather general conditions, the assumption $\sup L<\infty $ alone implies the attainment of a %maximum for $L$. The  in 
Proposition \ref{suplim}    develops  an idea from  L.R. Mead and N. Papanicolaou \cite{MP}, that we generalize to  our present  context.

\begin{proposition}\label{suplim}
Let $T$, $I$, $g$ and $\rho, a$ satisfy the hypotheses of Theorem {\rm \ref{lll}}.
Suppose also that $a (t)\! =\! \sum_{i\in I}c_i t^i$ and $\sum_{i\in I}c_i g_i \! >\! 0$. %Let $g=(g_i )_{i\in I}$ be arbitrary, and $L=L_{\rho ,a, g}$ be the corresponding Lagrangian.
If $\sup L_{\rho, a,g}\! <\! \infty$, 
then there is a   $\lambda^*$ on which the supremum is attained, $\sup L_{\rho ,a,g}\! =\! L_{\rho ,a,g}(\lambda^* )$.
\end{proposition}

{\it Proof}. Since $a$ is independent of $(t^i )_{i\in I\setminus \{ 0\}}$, $c_0 \not =0$.
 Set $c_{i0}=c_i$, $c_{ij}=\delta_{ij}$ ($i\in I$, $j\in I\setminus \{ 0\}$).
%The matrix $[c_{ij}]_{i,j\in I}$ is then invertible. 
A change of variables $\lambda \mapsto \tilde{\lambda}$: $\lambda_i =\sum_{j\in I}c_{ij}\tilde{\lambda}_j$ gives $L(\lambda )=\tilde{L} (\tilde{\lambda})\, :=\sum_{j\in I}\tilde{g}_j \tilde{\lambda}_j -\int_T e^{\sum_{j\in I}\tilde{\lambda}_j \tilde{a}_j -1}\rho \, dt$ where $\tilde{g}_j =\sum_{i\in I} c_{ij} g_i$ and $\tilde{a}_j =\tau_j /a$ for $\tau_j (t)=\sum_{i\in I}c_{ij}t^i$.
Then $\sup \tilde{L} =\sup L$. We prove  the attainment  for $\tilde{L}$. 
Denote $\tilde{\lambda}$,  $\tilde{a}_j$, $\tilde{g}$, $\tilde{L}$ by $\lambda$, $a_j$, $g$, $L$, respectively.
%Now $a_j (t)=\tau_j (t)/a(t)$ with $\tau_j (t)=\sum_{i\in I}c_{ij}t^i$ and 
Now  $a_{0}\equiv 1$, $g_0 >0$ and $(\tau_i )_{i\in I}$ are linearly independent on any subset of positive measure. Let $\mu =\rho dt$.  Since $\sup L<\infty$, $\mu (T)>0$.
%, $\mu (T)>0$.
 Set $\lambda =(\lambda_0 , \lambda ')$ where $\lambda '=(\lambda_i )_{i\in I\setminus \{ 0\}}$.
 Maximizing $L$ with respect to $\lambda_0$  gives $\alpha (\lambda '):=-\ln \int_T e^{\sum_{i\in I\setminus \{ 0\}}\lambda_i  a_i (t) -1}d\mu (t)$ such that  $\max_{\lambda_0 }L(\lambda_0 ,\lambda ')=L(\alpha (\lambda '),\lambda ')$. Consider the (convex) potential $f(\lambda '):=\int_T e^{\sum_{i\in I\setminus \{ 0\}}\lambda_i (a_i (t)-g_i )}d\mu (t)$  so that $\sup L<\infty $ $\Leftrightarrow$ $\inf f>0$. If $\inf f$ is attained at some $\lambda_* '$,
$\sup L$ will be attained at $(\alpha (\lambda_* '),\lambda_* ')$.  
 By (\ref{mon}), 
$
|\sum_{i\in I\setminus \{ 0\}}\lambda_i (a_i (t)-g_i )|\leq \| \lambda '\| ( 
 c\frac{\| t\|^m +1}{a(t)} + \| g\| )
$ where  $\| \lambda '\| \!  =\! \sum_{i\in I\setminus \{ 0\}}|\lambda_i ' |$, $\| g\| =\mbox{max}_{i\in I} |g_i |$ and $c$ is a constant. Then
 for every  sequence $\lambda_k ' \! =\! (\lambda_{ki} )_{i\in I\setminus \{ 0\}}$ such that  $\lim_{k}\lambda_k ' = \lambda '$  we have 
 %there is a constant $ct<\infty$ such that
$
e^{\sum_{i\in I\setminus \{ 0\}}\lambda_{ki} (a_i (t)-g_i )}\leq e^{\sup_{k}\! \| \lambda_k '\| ( 
 c\frac{\| t\|^m +1}{a(t)} \, + \| g\| )}
\in L^1 (T,\mu )
$. 
By Lebesgue's dominated convergence theorem,  $\lim_{k }f(\lambda_k ')=f(\lambda ')$.
Thus $f$ is continuous.

There is no  $\lambda '\not \! =\! 0$ such that  $p_{\lambda '}(t)\! :=\! \sum_{i\in I\setminus \{ 0\}}\lambda_i (\tau_i (t) /a(t) -g_i )\! \leq\! 0$ a.e. on $ T$,
for otherwise on the subset $Z : p_{\lambda '}(t)=0$ of $T$ we have $a(t)\sum_{i\in I\setminus \{ 0\}}\lambda_i g_i =\sum_{i\in I\setminus \{ 0\}}\lambda_i \tau_i (t)$; if $\sum_{i\in I\setminus \{ 0\}}\lambda_i g_i =0$, we get $\mu (Z)=0$ due to $\lambda '\not =0$; if
$\sum_{i\in I\setminus \{ 0\}}\lambda_i g_i \not =0$, we get again $\mu (Z)=0$ since $a$ is independent of  $(\tau_i )_{i\in I\setminus \{ 0\}}$ ($=(t^i )_{i\in I\setminus \{ 0\}}$). Hence $Z$ is negligible. Then on $T\setminus Z$,   $p_{\lambda '}(t)<0$, $e^{rp_{\lambda '}(t)}\leq 1$ ($r\geq 0$) and by Lebesgue's  theorem 
$f(r_k \lambda ')=\int_T e^{r_k p_{\lambda '} (t)}d\mu (t)\to 0$  
as  $r_k \to \infty$,
which is 
impossible since $\inf f>0$. %We have  used here that $a<\infty $ a.e.
 Then for any $\lambda '\not =0$
there are a constant $\delta =\delta_{\lambda ' }>0$ and  measurable subset $T_{\lambda '}\subset T$ with
$\mu (T_{\lambda '})>0$ such that $p_{\lambda '}(t)\geq \delta$ for all $t\in T_{\lambda '}$.
Hence $f(r\lambda ')\geq \int_{T_{\lambda '}}e^{rp_{\lambda '}(t)}d\mu (t)\geq e^{r\delta}\mu (T_{\lambda '})$. Then for every $\lambda '\not =0$,
$\lim_{r\to \infty}f(r\lambda ')=\infty$.  

There is a compact  $K\subset\mathbb{R}^{N-1}$ with  
$\inf f=\inf_{K}f$, for otherwise we could find a sequence of unit vectors $\lambda_k '$, and $r_k \to \infty$ 
such that  $\lim_{k\to \infty}f(r_k \lambda_k ')=\inf f$; 
 we can also assume there is a unit vector $\lambda '$ such that  $\lambda_k '\to \lambda '$. 
 Given  $r>0$,  $r\lambda_k '=s' \lambda_k '+(1-s')r_k \lambda_k '$ for $s'=\frac{r_k -r}{r_k -1}$ ($\to 1$ as $k\to \infty$) whence
 $  f(r\lambda_k ')\leq s' f(\lambda_k ')+(1-s')f(r_k \lambda_k ')$. Since $\sup_k |f(r_k \lambda_k ')|<\infty $ and $f$ is continuous, letting $k\to \infty$ we get $f(r\lambda ')\leq f(\lambda ')$ which is impossible because $\lim_{r\to \infty} f(r\lambda ')=\infty$.
%Thus $\inf f=\inf_{K}f$ for a compact $K$, and so
Since $\inf f$ is attained on $ K$, $\sup L$ will be attained. $\Box$
\vspace{3 mm}

% at $\lambda^* :=(\psi (\lambda_* '),\lambda_* ')$. 

A more explicit outcome  of  Theorem \ref{lll} and Proposition \ref{suplim} is  the Corollary \ref{c2} from below, that for small $\epsilon  $ is an approximate entropy maximization result.

\begin{corollary}\label{c2}
Let  $T\subset \mathbb{R}^n$ be a closed subset. Let $I\subset \mathbb{Z}_{+}^n $ be finite with $0\in I$. Fix $k\in \mathbb{Z}_+$ such that $\max_{i\in I}|i|< 2k+2$. Let $(g_i )_{i\in I}$ be a set of numbers with $g_0 =1$.
Fix also an arbitrary constant $\epsilon  >0$.
The following
statements {\em (a), (b), (c)} are equivalent:

{\em (a)} There exist functions $f\in L^{1}_{+} (T,dt)$ such that $\int_T |t^i |f(t)dt<\infty$ and
\begin{equation}\label{problemaa}
\int_T t^i f(t)\, dt =g_i \, ,\mbox{ }\mbox{ }\, i\in I;
\end{equation}

{\em (b)} There exists a particular solution $f_*$ of  {\em (\ref{problemaa})} maximizing the functional 
%the 
%functional $H:L_{+}^1 (T,dt)\to [-\infty ,\infty )$ given by 
$$H (f)=H_\epsilon (f)=-\int_T f\ln f\, dt -\epsilon   \int_T \| t\|^{2k+2} f\, dt;$$
%\frac{f}{\rho}\ln \frac{f}{\rho}\, \rho \, dt$$
%amongst all solutions;

{\em (c)} The associated Lagrangian $L=L_\epsilon$ from below satisfies $\, \sup L <\infty $
%=L (\lambda^{*}  )$ for a $\lambda^{*} =(\lambda_{i }^* )_{|i|\leq 2k}$ 
$$L (\lambda )=\sum_{i\in I}g_i \lambda_i -\int_T e^{\sum_{i\in I}\lambda_i t^i  -\epsilon  
\| t\|^{2k+2} } dt \, +1.$$ 
%is bounded from above, and attains its supremum in a point $\lambda^* $.

In this case: $\sup L$ is attained in a point $\lambda^*=(\lambda_{i}^* )_{i\in I}$,  both $f_*$ and $\lambda^*$  are uniquely determined,
$-H(f_* )=L(\lambda^* )$ and 
$$f_* (t)= e^{\sum_{i\in I}\lambda_{i}^* t^i -\epsilon   \| t\|^{2k+2}},$$
%resp. $f_{*\varepsilon} (t)= e^{\sum_{|i|\leq 2k}\lambda_{i}^* t^i -1-\epsilon \| t\|^{2k+2}}$)
in particular
$$ 
\int_T \, t^i \, e^{\sum_{j\in I}\lambda_{j}^* t^j }\, e^{-\epsilon   \| t\|^{2k+2}}  \, dt =g_i \mbox{ }\mbox{ }\mbox{ }\, (i\in I).
$$
%resp.
% $$ 
%\int_T \, t^i \, e^{\sum_{|j|\leq 2k}\lambda_{j}^* t^j -1-\epsilon \| t\|^{2k+2}}  \, dt =g_i \mbox{ }\mbox{ }\mbox{ }\mbox{ }\, (|i|\leq 2k).
%$$
\end{corollary}

{\it Proof}. Use Theorem \ref{lll} for $a(t)\equiv 1$ and  $\rho (t)=e^{-\epsilon   \| t\|^{2k+2} }$ ($t\in T$), which provides a Lagrangian $L_{\rho ,a,g}$ and point $\lambda^{*}_{\rho ,a,g}$ related to the present ones $L$, $\lambda^*$ by $L_{\rho ,a,g}(\lambda )=L(\lambda -\lambda^0 )$ and $\lambda^* =\lambda_{\rho ,a,g}^* -\lambda^0$ where $\lambda^0 =(\lambda_{i}^0 )_{i\in I}$ with $\lambda_{i}^0 =\delta_{i0}$. 
Then Proposition \ref{suplim}  applies, since $\sum_{i\in I}c_i g_i =g_0 >0$. $\Box$
\vspace{3 mm}

\noindent {\bf Remarks} Considering   perturbations  $H_\epsilon $  of the entropy $H_0$ as above  might automatically provide enough control in the tails of a maximizing sequence to guarantee convergence by known arguments \cite{rod}, \cite{Moreau}.  More specifically, a maximizing sequence converges in weak $L^1$, and if the dominant moments of order $2k+2$ were  bounded,  one could show that  moments of lower order $\leq 2k$ converge.  For this argument the author is indebted to the referee, who  legitimately suggested that   Corollary \ref{c2} and Theorem \ref{lll} may  have  shorter proofs by this standard method, and moreover   results like Corollary \ref{c2} are rather known 
\cite{gmc}, see below. This seems to be true indeed, because maximizing $H_0 -\epsilon \int_T \| t\|^{2k+2}f(t)dt$ should provide a certain brake on the growth of the moments of order $2k+2$. 
However, within the quite technical hypotheses of Theorem \ref{lll} we could not find an obvious argument to get apriori bounds on the $2k+2$ moments, and  for the sake of completeness 
%.
%This is rather credible since
%successive integration by parts in equations like $\int \frac{\partial \, }{\partial t_j}(t_j e^{\sum_{|i|\leq 2k}\lambda_i t^i -\| t\|^{2k+2}})dt=0$ can provide recurrence relations involving the %moments, that lead to estimates of the form $|g_j |\leq ct_k \, \sqrt{1+\| \lambda \|^2 }\sqrt{1+\| g\|^2 }$, $|j|=2k+2$. Since the variable factor $ \sqrt{1+\| \lambda \|^2 }$ raises %complications to this aim, and for 
 we kept  our initial proofs.      Theorem \ref{lll}  can  at least  unify and cover also  other known cases, like  $T:=\, $compact, 
$\epsilon =0$ \cite{Leww}, see also \cite{io} (setting $a(t),\, \rho (t)\, \equiv 1$ on $T$) or, to some extent, 
[Theorem 8, \cite{bla}] setting
$a(t)=(\| t\|^2 +1)^k$, $\rho (t)=c\| t\|^{2-n}(\| t\|^2 +1)^{-3/2}$ ($n\! \geq \! 2$); we omit the details. Another  application is for example Corollary \ref{cc} from below.

In principle,
one could numerically maximize such Lagrangian functions $L$  to  obtain a vector $\lambda^*$ and so a density $f_*$. Solving such dual problems (usually by  Newton's method) turns to be the basic technique to this aim. The main effort is then  to deal with the computational cost of approximating  multiple integrals (like $\int_T t^j e^{\sum_{i\in I}\lambda_i t^i }\rho (t)dt$ if $a\equiv 1$, for instance)  needed for the gradient of $L$
 \cite{junk}, \cite{Levermore},   \cite{lassrec}, \cite{lassrec2},   \cite{bdl},    \cite{abr}, \cite{gmc}. 
%\vspace{3 mm} 

Versions of Corollary \ref{c2} were indeed known, for instance in C.P. Groth, J.G.  McDonald \cite{gmc}  that used an ansatz like $f_* (t)=e^{\sum_{i\in I}\lambda_i t^i}\rho (t)$  to derive a moment closure in the context of kinetic equations, when $t\in \mathbb{R}^3$ stands for the gas velocity \cite{cerc}. Their density  $\rho (t)=e^{\sigma (t)}$ with a certain negative $\sigma (t)$ of order $|\sigma (t)|\geq ct.\,   \| t\|^{2k+2}$, see [eqs. (51), (52), \cite{gmc}]  can be viewed as a "window" function that attenuates the distribution at high velocities.  This type of modification to the maximum-entropy moment distribution has been  proposed by J.D. Au \cite{au}, and M. Junk \cite{junku} and   prevents the existence of very small packets of very  fast particles that, as mentioned in   \cite{junke},  are the basic reason for non-solvability of the maximum entropy problem associated with Euler equations \cite{cerc}.

Our proposed technique has also  similarities with the method used in  A. Tagliani, \cite{tag} that deals with the particular case of the Hamburger moments problem for $T=\mathbb{R}$.
% Maximum entropy in the Hamburger moments problem, J. Math. Phys., 35, 5087, 1994 by A. Tagliani....................
\vspace{2 mm}

%   provided the corresponding integrals can be well approximated.

\begin{corollary}
\label{cc}Let $T\subset \mathbb{R}^n$ be  closed,  
$k\in \mathbb{Z}_+$, $I=\{ i\in \mathbb{Z}_{+}^n :|i|\leq 2k\} $ and $(g_i )_{|i|\leq 2k}$  a set of reals  with $g_0 =1$. 
The statements {\em (a), (b)}  are equivalent:

{\em (a)} There exists an  $f\in L^{1}_{+} (T,dt)$ such that
$$
\int_T t^i f(t)\, dt =g_i \mbox{ }\mbox{ }\mbox{ }\, (|i|\leq 2k);
$$
%\end{equation}

{\em (b)}  $L(\lambda ):=\sum_{|i|\leq 2k}g_i \lambda_i -\int_T e^{\frac{\sum_{|i|\leq 2k}\lambda_i t^i }{\| t\|^{2k} +1}}e^{-\| t\|^2 -1} dt$ is bounded from above.  

\noindent In this case, $L$ attains its maximum in a point $\lambda^* =(\lambda_{i}^* )_{|i|\leq 2k}$ and
 $$f_* (t):= \frac{1}{\| t\|^{2k}+1 }\, e^{\sum_{|i|\leq 2k}\lambda_{i}^* \frac{t^i }{\| t\|^{2k}+1}}\, e^{-\| t\|^2 -1} $$ 
satisfies $\int_T t^i f_* (t)\, dt =g_i $ ($|i|\leq 2k$).
\end{corollary}

{\it Proof}.   Use Theorem \ref{lll} and  Proposition \ref{suplim} for $a=\| t\|^{2k} +1$, $\rho =e^{-\| t\|^2}$. $\Box$
\vspace{3 mm}

%\newpage
% We consider  Lagrangians $L_{\epsilon} $ ($\epsilon \geq 0$), for  $H\approx -\int_{}f\ln f\, dt $ and $\, I\! :\! |i|\!\leq\! 2k\, $. 

\noindent  {\bf Notation} For  $g\! =\! (g_i )_{|i|\leq 2k} $ having representing densities on $\mathbb{R}^n$,  let 
$\lambda^* \! =\! \lambda^{*}_g \!  =\! (\lambda_{i}^* )_{|i|\leq 2k}$ denote the  vector maximizing  $L_0 (\lambda )\! =\! \sum_{|i|\leq 2k}g_i \lambda_i \! -\! \int e^{\sum_{|i|\leq 2k}\lambda_i t^i} dt$. Set $p_g (t)\! =\! \sum_{|i|\leq 2k}\lambda_{i}^* t^i$. Then $\sum_{|i|=2k}\lambda_{i}^* t^i \! \leq \! 0$ for all $t\! \in \!\mathbb{R}^n$ (use  $\int_{\mathbb{R}^n }\! e^{p_g}dt\! <\! \infty$ and polar coordinates).
%\vspace{1 mm}
 Let $G$ be the set of all such $g$, with the  property $\sum_{|i|=2k}\lambda_{i}^* t^i \! <\! 0$ for 
 $t\! \not =\! 0$. Then (see  \cite{Levermore})
$G$ is dense and open in the set of all $g$ having representing densities, consists of data $g$ for which $\lambda^* $  
does provide a representing density $f_* \! =\! e^{p_g}$ of $g$ maximizing $H(f)\! =\! -\int f\ln f dt$,  
and the map $G\ni g\mapsto \lambda_{g}^*$ is $C^\infty$-diffeomorphic. Let  $n,\, k=2$, whence $\mbox{card}\, \{ i\! :\! |i|\! \leq \! 2k\}\! =\! 15$.
Let $x=(x_i )_{i\in \mathbb{Z}_{+}^2 ,\, |i|\leq 4}$ denote the variable in $\mathbb{R}^{15}$.
%If $g\in G$,  both $\lambda^*_{(4,0)}$, $\lambda^*_{(0,4)}\, <0 $.
Let $G_0 \! =\! \{  g\! \in\!  G\! :  \det A \lambda^{*} \! \not =\! 0\}$ where
 $A=Ax$ is the  matrix in  (\ref{calcul}). Then $G_0$ is dense and open in $G$. Given $g\in G$, we may set $g_j :=\int t^j e^{p_g (t)}dt$ for $|j|\geq 5$.
\vspace{3 mm}

Proposition \ref{p1} from below   is reminiscent to J.B. Lasserre [Lemma 2, \cite{lassrec}], see also \cite{lassrec} or [Lemma 2, \cite{mpb}], where  similar recurrences  were obtained.  The idea in the case $n=1$ is to compute integrals like $\int \frac{d\, }{dt}(tf_* (t))dt$ via Leibniz-Newton's formula. In our  case $n=2$ the basic idea is the same, but a careful application of Stokes' theorem will be required in the proof. Although such  calculation is not a practical method itself for determining if a moment set comes from an underlying density, it could help to the approximation of $\lambda^*$
when used along with suitable numerical techniques - see for example  the use of Newton's method as in \cite{lassrec} together with the semidefinite programming methods for gradient and  Hessian  computation from \cite{bdl}, \cite{lassrec2}.

%- for instance, Newton's method to maximize the Lagrangian $L$  \cite{lassrec}, and  \cite{bdl}. 

\begin{proposition}\label{p1} Let $n, k\! =\! 2$ and $g\! \in\! G_0$. The higher order moments $(g_j )_{|j|\geq 5}$ of the maximum entropy density $e^{p_{g}}$ can be expressed by relations of the form
$$g_j =\sum_{|i |\leq 4}r_{j i} (\lambda^* )\, g_i \mbox{ }\mbox{ }\mbox{ }\mbox{ }\,
\mbox{ }(j\in \mathbb{Z}^{2}_+ , \, |j|\geq 5)$$  where  $r_{ji}=r_{ji}(x)$ are  universal rational functions, see {\em (\ref{calcul}) -- (\ref{ratz})}.
\end{proposition}

{\it Proof}.
It suffices to prove that for any $j_0 \in \mathbb{Z}_{+}^2$ with $|j_0 |\geq 5$ there are rational functions $c_{j_0 i}=c_{j_0 i}(\lambda^* )$, for $|i|<|j_0 |$, such that
$g_{j_0} =\sum_{|i|<|j_0 |}c_{j_0 i}g_i$ and then proceed inductively. Set $|j_0 |=l+1$ for $l\geq 4$ and 
denote $\lambda^* =(\lambda^{*}_i )_{|i|\leq 4}$ by $x=(x_i )_{|i|\leq 4}$.  Set $x_\kappa =0$ if $\kappa \not \geq 0$. 
Let $p=p_g$, namely $p(t)=\sum_{|\iota |\leq 4}x_\iota t^\iota$. We will find a polynomial 
$\pi (t)=\sum_{|i|\leq l}c_{j_0 i}t^i$ and a differential 1--form $\omega =e^p (udt_1 +vdt_2 )$ with   $u,v$ polynomials, depending on $j_0$, such that  $d\omega (t)=(t^{j_0} -\pi (t))e^{p (t)}dt_1 \wedge dt_2$. By Stokes' theorem on  disks $D_r$ of center 0 and radius $r$, $\int_{D_r}d\omega =\int_{\partial D_r}\omega \to 0$ as $r\to \infty$ since $ue^p ,ve^p $ are rapidly decreasing ($g\in G$). Hence
$\int_{\mathbb{R}^2}(t^{j_0} -\pi (t))e^{p (t)}dt_1  dt_2 =0$ which is the desired conclusion. The condition on $\omega$ means that $L=L(u,v):=v\partial_1 p -u\partial_2 p +\partial_1 v -\partial_2 u $ where $\partial_m =\partial /\partial t_{m}$ ($m=1,2$) satisfies $L =t^{j_0} -\pi $. We let  $u(t)=\sum_{|j|=l-2}a_j t^j$, $v(t)=\sum_{|j|=l-2}b_j t^j$ with $a_j \! =\! a_j (x), b_j \! =\! b_j (x)$
rational functions to be determined.  Set $e_1 =(1,0), e_2 =(0,1)$. In  degree $l+1$,
the equation $L=t^{j_0} -\pi $ gives
$
\sum_{|j|=l-2,\, |\iota |=4,\,  \iota_1 \geq 1}b_j \iota_1 x_\iota t^{j+\iota -e_1}-
\sum_{|j|=l-2,\, |\iota |=4,\,  \iota_2 \geq 1}a_j \iota_2 x_\iota t^{j+\iota -e_2}=t^{j_0}
$.
Change the summation  indices by $i=j+\iota -e_{1,2}$ and identify
the coefficients of $t^i$ with $i\geq 0,\, |i|=l+1$. Then
\begin{equation}\label{myb}
\sum_{|\iota |= 4,\, ( e_1 \leq \iota \leq i+e_1 )}\! \! \iota_1 x_\iota b_{i+e_1 -\iota}-
\! \! \sum_{|\iota |= 4,\,  (e_2 \leq   \iota \leq i+e_2 )}\! \! \iota_2 x_\iota a_{i+e_2 -\iota}=\delta_{i j_0 }\mbox{ }
\mbox{ }\mbox{ (}|i|=l+1\mbox{)}\end{equation}
% for
%all $i\geq 0$ with $|i|=l+1$,
where $\delta_{i j_0}$ is Kronecker's symbol. 
The summation conditions in the brackets from above may  be omitted, since the  terms outside the respective ranges vanish formally  due to either  $\iota_{1,2}=0$, or
 $a_j ,b_j ,x_\kappa =0$ whenever $j,\kappa \not \geq 0$. Once we have such $u,v$,  $\pi$ is determined from $L=t^{j_0}-\pi$ by gathering all terms of degree $\leq l$ in $-L$.
% (we look for $u,v$ polynomial).
We solve (\ref{myb})  in the Appendix, that provides also an algorithm for computing $c_{j_0 i}$, $r_{ji}$
via 
the formulas (\ref{calcul}) -- (\ref{ratz}).
 $\Box$
\vspace{3 mm}

Corollary \ref{exempluu} below is an attempt to solve a  maximum entropy problem  by means of a system of ordinary differential equations (\ref{sistem}) without computing multiple integrals.
 However this is rather a theoretic result, since it requires an accurate solution of (\ref{sistem}) - that is, small increments $\Delta s$ and so, a large number $1/\Delta s$ of iterations.

%To these aims, the approximation
%f the integrals like $\int t^i f_* (t)dt$ necessary to compute the entries of the matrix $H$ became possbile by now for $n\geq 3$ and small values of $k$ by various numerical techniques, see %for instance \cite{lassrec}, \cite{bdl}.

% for the
% functions $\gamma_j \! =\! \gamma_j (x)$  defined in (\ref{nrr}), see also (\ref{ratz}).

\begin{corollary}\label{exempluu}
 Let $n,k\! =\! 2$ and $g, g_0 \! \in \!\mathbb{R}^{15}$  such that $sg  + (1\! -\! s) g_0 \in G_0$ for all $s\! \in\!  [0,1]$, where
 $g_0$ has a known $\lambda^{ *}_{g_0}$. Set $\Gamma_i (x,s)=sg_i+(1-s)(g_0 )_i$ for $|i|\leq 4$ and 
 $\Gamma_j (x,s)=  \sum_{|i |\leq 4}r_{j i} (x)\, (sg_i +(1-s)(g_0 )_i )$ for  $|j|\geq 5$ where $x=(x_i )_{|i|\leq 4}\in \mathbb{R}^{15}$. 
 The system of  ordinary differential equations
 \begin{equation}\label{sistem}
 \sum_{|j|\leq 4}\Gamma_{i+j}\, (x(s),s)\, \frac{dx_j}{ds\, {\, }\, }\, (s)=g_i -(g_{0})_i \mbox{ }\mbox{ }\mbox{ }(|i|\leq 4);
 \mbox{ }\mbox{ }\mbox{ }\, x(0)\! =\! \lambda^{*}_{g_0}
 \end{equation}
 has a $C^\infty $ solution $x\! =\! x(s)$,  defined  on a neighborhood of $[0,1]$,  the matrix $[\Gamma_{i+j}(x(s),s)]_{|i|,|j|\leq 4}$ is  defined and invertible for all $s\in [0,1]$, and  we have $x(1)=\lambda^{*}_g $.
 \end{corollary}
 
 {\it Proof}.   Since $G_0$ is open,  the point
 $g(s):=s g +(1-s)g_{0} $ is in $ G_0$  (in particular, has representing densities) for every  $s$ in a neighborhood of $[0,1]$. Set $g(s)=(g_i (s))_{|i|\leq 4}$. Since $g(s)\in G$, %($\subset G$), 
 it has a  $\lambda^{*} =\lambda^{*}_{g(s)}$ maximizing $L_{0 ,\, g(s)}$.
 Let $x(s)=\lambda^{*}_{g(s)}$. Write $x(s)=(x_\iota (s))_{|\iota |\leq 4}$. Then $p_{g(s)} (t)=\sum_{|\iota |\leq 4}x_\iota (s)t^\iota$. The $H$--maximization holds and $e^{p_{g(s)}}$ is a representing density for $g(s)$,
 \begin{equation}\label{ggg}g_i (s)= \int_{\mathbb{R}^2} t^i e^{\sum_{|\iota |\leq 4}x_\iota (s)t^\iota }dt \mbox{ }\mbox{ } \mbox{ (}|i|\leq 4\mbox{)}.\end{equation}
 Denote by $g(s)_j$ for $|j|\geq 5$ the moments of higher order of $e^{p_{g(s)}}$, namely $g(s)_j :=\int t^j e^{p_{g(s)} (t)}dt$ ($|j|\geq 5$).
 Since  the map $G\ni \tilde{g}\mapsto \lambda^{*}_{\tilde{g}}$ is diffeomorphic, $x(\, \cdot \, )$ is smooth and so we may apply  $d/ds$ to the equalities (\ref{ggg}), whence 
 $$
 g_i -(g_0 )_i = \sum_{|j|\leq 4}\int t^{i+j } e^{\sum_{|\iota |\leq 4}x_\iota (s)t^\iota }\, \frac{dx_j }{ds}(s)= \sum_{|j|\leq 4}g(s)_{i+j}\frac{dx_j }{ds}(s); $$ of course $g(s)_i =g_i (s)$ if $|i|\leq 4$. Also %for all $j\in \mathbb{Z}_{+}^2$ with $|j|\geq 5$ we have
%$g(s)_j =\gamma_j (x)$  for $x=\lambda^{*}_{g(s)}$. 
%Thus $g(s)_j$  can be  computed by plugging $\lambda^{*}_{g(s)}$ ($=x(s)$) and $g(s)$ in the right hand side of. 
 $g(s)_{j} =\sum_{|i|\leq 4}r_{ji}(x(s))g_i (s) =\Gamma_{j} (x(s),s)$ ($|j|\geq 5$) by Proposition \ref{p1}. Then 
we obtain the differential equations (\ref{sistem}) on a neighborhood of $[0,1]$.  
The denominators of  $r_{ji}(x)$ do not vanish on the set $\{ x(s):0\! \leq\!  s\! \leq\!  1\} $ and so  $\Gamma_{i+j}(x(s),s)$ are defined for $0\! \leq\! s\!\leq\! 1$.
Each matrix $[\Gamma_{i+j}(x(s),s)]_{|i|,|j|\leq 4} = [\int t^{i+j}e^{p_{g(s)}(t)}dt]_{|i|,|j|\leq 4}$ is positive definite and so invertible.
 By (\ref{ggg}), $g_i =\int t^i e^{\sum_{|\iota |\leq 4}x_\iota (1)t^\iota }dt$ ($|i|\leq 4$). Due to the uniqueness of the critical point  of  the Lagrangian $L_{0,g}$ we derive $x(1)=\lambda^{*}_{g}$. $\Box$
\vspace{3 mm}

\noindent {\bf Remark 10}  Since  $\Gamma_{j} (x(s),s)=g(s)_{j} $  for $|j|\geq 5$ where $g(s)=sg+(1-s)g_0$, all the entries of the matrix $\Gamma =[\Gamma_{i+j}(x(s),s)]_{|i|,|j|\leq 4}$ of the system (\ref{sistem}): $\Gamma (x(s),s)\cdot \frac{dx}{ds}(s)=g-g_0$ are moments and can be computed inductively by  linear recurrences  $g(s)_{j_0} =\sum_{|i|\leq l}c_{j_0 i}(x(s))g(s)_i$ ($|j_0 |=l+1$), see  (\ref{nrr}), using the concrete formulas (\ref{calcul}), (\ref{coeff}) of  $c_{j_o i}$; the explicit
 formulas of $r_{ji}$  from (\ref{ratz}) are not needed to this aim. Moreover, for each $l$ the calculations of $g_{j_0}$ ($|j_0 |=l+1$) are independent of each other. We may consider any $g_0 \in G_0$, for instance the set of moments up to the 4th order of $e^{-t_{1}^4 -t_{2}^4}$.  Also fast inversion algorithms exist for such Hankel matrices  $\Gamma$. 
 Then for problems of reasonable size  
one can use   numerical methods for systems of ordinary differential equations   to obtain $\lambda^{*}_g$ ($=x(1)$).
 \vspace{2 mm}

The author is indebted to one of the referees for  Remark 11 from below. 
\vspace{2 mm}

\noindent {\bf Remark 11} Due to the accuracy needed for its solution the system (\ref{sistem})   is not a very practical way of computing $\lambda^*$, comparing to the  more efficient Newton's method or its many variants  \cite{lassrec}, \cite{lassrec2}, \cite{junk}, \cite{Levermore}, \cite{bdl} for the dual problem of maximizing $L$. Actually, 
one could use (\ref{sistem}), written in the form $dx/ds =  H(x(s)) (g - g_0 )$, to iteratively find $x(1)$ as follows.  Let the increment $ \Delta s = 1$.  Then a forward Euler solve of the ODE gives
$
 x_1 = x(0) - H(x(0)) (g - g_0 )
$.
 The value $x_1$ is an estimate of $x(1)$, but not exactly, so one can repeat the process with $x_1$ as the new initial condition and a new moment $g_1$ which is computed from $x_1$.  Doing this $k$ times gives
$
 x_k = x_{k-1} - H(x_k ) (g - g_k )
$.
 If we identify $H$ with the Hessian of the dual problem, then this is just Newton's method. When  is well conditioned, a standard Newton method for the dual problem will take only a handful of such iterations.
\vspace{3 mm}
 
\noindent {\bf Appendix. The functions ${\bf r}_{\bf{j\, i}}$ } 
We give an algorithm to recurrently compute $r_{ji},c_{ji}$, in particular solve (\ref{myb}) to finish the proof of Proposition \ref{p1}.
%Letting $i=(l+1-k,k)$ for $0\leq k\leq l+1$ in (\ref{myb}) gives $l+2$ equations. 
Set $\delta_k =\delta_{(l+1-k,k)\, j_0 }$ for $0\leq k\leq l+1$. 
Let $\alpha_k =a_{(l-2-k,k)}$, $\beta_k =b_{(l-2-k,k)}$ for $0\! \leq \! k\! \leq\! l\! -\! 2$. Thus 
 $\alpha_k ,\beta_k =0$ for  $k<0$,   $k\geq l-1$. Also $x_\kappa =0$ if $\kappa \not \geq 0$.
Change the summation indices  in (\ref{myb}) by $j\! =\! i\! +\! e_{1,2}\! -\! \iota$ ($ \geq\! 0$). Then (\ref{myb}) becomes
$\sum_{|j|=l-2, (j_2 \leq i_2 ) }(i_1 -j_1 +1)x_{i-j+e_1}b_j -\sum_{|j|=l-2,( j_2 \leq i_2 )}(i_2 -j_2 +1)x_{i-j+e_2}a_j =\delta_{ij_0}$ 
where the (redundant) condition $j_2\leq i_2$ follows from $j\leq i$, that comes from $\iota \geq e_{1,2}$.  For every $i=(l+1-k,k)$ with 
$0\leq k\leq l+1$, we have the equivalence $(j\geq 0,|j|=l-2, j_2 \leq i_2 )\, \Leftrightarrow   \, j=(l-2-p,p)$ for $0\leq p\leq k$ and  
hence the $l+1$ equations in (\ref{myb}) become now, respectively,
\begin{equation}\label{cumva}
\sum_{p=0}^k [ (4\! +\! p\! -\! k)x_{(4+p-k,k-p)}\beta_p \! - \! (k\! -\! p\! +\! 1)x_{(3+p-k,k-p+1)}\alpha_p ] \! =\! \delta_k ,\, 0\! \leq\! k\! \leq\!  l+1.
\end{equation}
If $l\geq 5$,  let $\alpha_0 ,\ldots ,\alpha_{l-5}=0$ and define  $\beta_0 ,\ldots ,\beta_{l-5}$ inductively by
$
4x_{40}\beta_k = -\sum_{p=0}^{k-1}(4+p-k)x_{(4+p-k,\, k-p)}\beta_p +\delta_k$ ($0\leq k\leq l-5$)
where $\sum_\emptyset :=0$. Note that $x_{40}<0$ since $g\in G$. This fulfills (\ref{cumva}) for $0\leq k\leq l-5$. 
%Once we have $\alpha_k ,\beta_k $  for $0\leq k\leq l-5$, 
Last 6 equations in (\ref{myb}) ($l-4\leq k\leq l+1$ in (\ref{cumva})) will provide  $\alpha_k ,\beta_k $ ($l-4\leq k\leq l-2$).  
If $l=4$,  skip this step and go directly to the linear $6\times 6$ system
(in this case (\ref{calcul})  for $y,z,w=0$). In any case, 
%the form (\ref{cumva}) of the equations  is not suitable to this aim. Instead, 
we let now
  $i=(l+1-k,k)$ for  $0\leq k\leq l+1$ in (\ref{myb}). We have $i+e_{1}-\iota =(l+2-k -\iota_1 ,k-\iota_2 )$ and $i+e_2 -\iota =(l+1-k-\iota_1 ,k-\iota_2 +1 )$. 
Last 6 equations in (\ref{myb})
 become  $
\sum_{|\iota |=4} \! \iota_1 x_\iota \beta_{k -\iota_2 }
\! -\!  \sum_{|\iota |=4} \! \iota_2 x_\iota \alpha_{k-\iota_2 +1}=\delta_k $  ($l-4\leq k\leq l+1$), see below.
The brackets $(\mbox{ })$ border  quantities
  already known in terms of $\beta_0 ,\ldots ,\beta_{l-5}$. The 
markers $\lfloor\mbox{ } \rceil$  border sums of terms that are null due to
 $\iota_{1,2}=0$,
$\alpha_k ,\beta_k =0$ ($k\geq l-1$) or $\alpha_k =0$ ($0\leq k\leq l-5$):
 %write it for $k=l-4$, then  copy the equation increasing the indices of $\alpha_k ,\beta_k $ by 1 each time: 
$$
\begin{array}{l}4x_{40}\beta_{l-4}+(3x_{31}\beta_{l-5}+2x_{22}\beta_{l-6}+1x_{13}\beta_{l-7}+0x_{04}\beta_{l-8}) \hspace{3 cm} \\
\hspace{1.2 cm} \lfloor -0x_{40}\alpha_{l-3}\rceil -1x_{31}\alpha_{l-4}\, \lfloor -2x_{22}\alpha_{l-5}-3x_{13}\alpha_{l-6}-4x_{04}\alpha_{l-7} 
\, \rceil =\delta_{l-4}\end{array} 
$$
$$
\begin{array}{l}4x_{40}\beta_{l-3}+3x_{31}\beta_{l-4}+(2x_{22}\beta_{l-5}+1x_{13}\beta_{l-6}+0x_{04}\beta_{l-7}) \hspace{3 cm} \\
\hspace{1.2 cm} \lfloor -0x_{40}\alpha_{l-2}\rceil -1x_{31}\alpha_{l-3}-2x_{22}\alpha_{l-4}\, \lfloor -3x_{13}\alpha_{l-5}-4x_{04}\alpha_{l-6}\, \rceil  =\delta_{l-3}\end{array} 
$$
$$
\begin{array}{l}4x_{40}\beta_{l-2}+3x_{31}\beta_{l-3}+2x_{22}\beta_{l-4}+(1x_{13}\beta_{l-5}+0x_{04}\beta_{l-6}) \hspace{3 cm} \\
\hspace{1 cm} \, \lfloor -0x_{40}\alpha_{l-1}\, \rceil -1x_{31}\alpha_{l-2}-2x_{22}\alpha_{l-3}-3x_{13}\alpha_{l-4}\, \lfloor -4x_{04}\alpha_{l-5}\, \rceil  
=\delta_{l-2}\end{array} 
$$
$$
\begin{array}{l}\, \lfloor 4x_{40}\beta_{l-1}+\, \rceil \, 3x_{31}\beta_{l-2}+2x_{22}\beta_{l-3}+1x_{13}\beta_{l-4}+\lfloor 0x_{04}\beta_{l-5}\rceil \hspace{3 cm} \\
\hspace{1.5 cm} \, \lfloor -0x_{40}\alpha_{l}-1x_{31}\alpha_{l-1}\, \rceil -2x_{22}\alpha_{l-2}-3x_{13}\alpha_{l-3}-4x_{04}\alpha_{l-4} 
=\delta_{l-1}\end{array} 
$$
$$
\begin{array}{l}\, \lfloor 4x_{40}\beta_{l}+3x_{31}\beta_{l-1}+\, \rceil \, 2x_{22}\beta_{l-2}+1x_{13}\beta_{l-3}+\lfloor 0x_{04}\beta_{l-4}\rceil  \hspace{3 cm} \\
\hspace{1.5 cm} \, \lfloor -0x_{40}\alpha_{l+1}-1x_{31}\alpha_{l}-2x_{22}\alpha_{l-1}\, \rceil -3x_{13}\alpha_{l-2}-4x_{04}\alpha_{l-3} 
=\delta_{l}\end{array} 
$$
$$
\begin{array}{l}\, \lfloor 4x_{40}\beta_{l+1}+3x_{31}\beta_{l}+2x_{22}\beta_{l-1}+\, \rceil \, 1x_{13}\beta_{l-2}+\lfloor 0x_{04}\beta_{l-3}\rceil  \hspace{3 cm} \\
\hspace{1.5 cm} \, \lfloor -0x_{40}\alpha_{l+2}-1x_{31}\alpha_{l+1}-2x_{22}\alpha_{l}-3x_{13}\alpha_{l-1}\, \rceil -4x_{04}\alpha_{l-2} 
=\delta_{l+1}.\end{array} 
$$
Set $
y\! =\! -3x_{31}\beta_{l-5} \! -\! 2x_{22}\beta_{l-6}\!  -\! x_{13}\beta_{l-7}$, $ z\! =\! -2x_{22}\beta_{l-5}\! -\! x_{13}\beta_{l-6}$ and $w\! =\! -x_{13}\beta_{l-5}
$.
 We easily read from above that $\alpha_k ,\beta_k$ for $k\! =\! l-4,\, l-3,\, l-2$  are given by

\begin{equation}\label{calcul}
\left[ \, \begin{array}{ccccccc} 4x_{40} & 0&0&& x_{31}&0&0\\ 
3x_{31}& 4x_{40}& 0&& 2x_{22}& x_{31}& 0\\
2x_{22}& 3x_{31}& 4x_{40}&& 3x_{13}& 2x_{22}& x_{31}\\
&&&&&&\\
x_{13}& 2x_{22}& 3x_{31}&& 4x_{04}& 3x_{13}& 2x_{22}\\
0& x_{13}& 2x_{22}&& 0& 4x_{04}& 3x_{13}\\
0& 0&x_{13}&&0&0&4x_{04}
\end{array}\, \right] \!   \left[ \, \begin{array}{l}\beta_{l-4}\\ \beta_{l-3}\\ \beta_{l-2}\\ \\ -\alpha_{l-4}\\ -\alpha_{l-3}\\
-\alpha_{l-2} \end{array} \, \right] \! \! =\! \!  
\left[  \begin{array}{l}y+\delta_{l-4}\\ z+ \delta_{l-3}\\ w+\delta_{l-2}\\ \\ \delta_{l-1}\\ \delta_{l}\\
\delta_{l+1} \end{array}\right] 
\end{equation}
(note also that $g\in G_0$). We have  $a_j, b_j$, and so $u, v$ such that $\mbox{deg}\, (L(u,v)-t^{j_0})\leq l$.  Now $\pi =t^{j_0}-L$ is determined by summing  the terms of degree $\leq l$ in $-L$. For $m=1,2$ set $K_m = \{ (j,\iota ): |j|=l-2, |\iota |\leq 3, \iota_m \geq 1\}$. Then 
$$
\pi \! =\! \! \! \! \! \sum_{(j,\, \iota )\in K_2 }\! a_j \iota_2 x_\iota t^{j+\iota -e_2}\! -\! 
        \! \! \!\! \sum_{(j,\, \iota )\in K_1 }\! b_j \iota_1 x_\iota t^{j+\iota -e_1}\! 
+\! \! \! \! \!\sum_{|j|=l-2, j_2 \geq 1}\! j_2 a_j t^{j-e_2}\! -\! 
 \! \! \! \!\sum_{|j|=l-2, j_1 \geq 1}\! j_1 b_j t^{j-e_1}.$$
 For any $i\geq 0$ with $|i|\leq l$, 
 the coefficient of $t^i$ in the sum $\Sigma_{K_2}$ from above is $\sum_{(j,\, \iota )\in K_2 (i) }a_j \iota_2 x_\iota$ 
 where $K_2 (i)= \{ (j,\iota )\in K_2 :j+\iota -e_2 =i\}$. The map $K_2 (i)\ni (j,\iota )\mapsto i-j$ is bijective onto 
  $I_i :=\{ \kappa \geq 0 :\kappa \leq i,|\kappa |=|i|+2-l\}$. Then we may use it to change the summation index by $\kappa =i-j$ and get the coefficient of $t^i$ in $\Sigma_{K_2}$  as $\sum_{\kappa \in I_i} (\kappa_2 +1)a_{i-\kappa }x_{\kappa +e_2}$. Similarly, the coefficient of $t^i$ in $\Sigma_{K_1}$ is $\sum_{\kappa \in I_i}(\kappa_1 +1)b_{i-\kappa}x_{\kappa +e_1}$.
% By changes of index $j=i+e_m-\iota$, $\kappa =\iota -e_m$ in the first two sums from above we get t
The coefficient
  $c_{j_0 i}$ ($=\, $a rational function  $c_{\, lj_0 i}(x)$ of $x$, actually)  of $t^i$ in $\pi (t)$ is then
\begin{equation}\label{coeff}
c_{j_0 i} =\sum_{\kappa \in I_i}[(\kappa_2 +1)x_{\kappa +e_2}a_{i-\kappa} -(\kappa_1 +1)x_{\kappa +e_1}b_{i-\kappa} ]+d_{j_0 i}\mbox{ }\mbox{ }\mbox{ }\mbox{ (}|i|\leq l\mbox{)}\end{equation}
where $
d_{j_0 i}\! =\! (i_2 +1)a_{i+e_2 }\! -\! (i_1 +1)b_{i+e_1 }$ if $ |i|\! =\! l-3$, and 0 otherwise.
%\mbox{)}.
%$
%For every $l\geq 4$ and $j_0$ with $|j_0 |=l+1$ we have the set of functions $c_{j_0 i}(x)$, $|i|\leq l$ from above.
%Therefore we found $c_{j_0 i}$ such that $g_{j_0}=\sum_{|i|\leq l}c_{j_0 i}g_i$ for all $l\geq 4$ and $j_{0}$ with $|j_0 |=l+1$. 
%We emphasize now the dependence $c_{j_0 i}=c_{j_0 i}(x)$, on the vector 
We have
\begin{equation}\label{nrr}
g_{j_0} =\sum_{|i|\leq l}c_{j_0 i}(x)\, g_i \mbox{ }\mbox{ }\mbox{ }\mbox{ }\mbox{ (} |j_0 |=l+1,\, l\geq 4\mbox{)}.\end{equation}
%Of course $ \gamma_j (x)=g_j $ when $x=\lambda^* $. %where both $j_0$  and $ l $ are now both variable. 
Successive compositions of the mapping $(g_i )_{|i|\leq l}\mapsto ((g_{j_0})_{|j_0 |=l+1},\, (g_i )_{|i|\leq l})=$ $(g_i )_{|i|\leq l+1}$ given by (\ref{nrr}) for $l=4,5,\ldots $ provide us with  $r_{ji}(x)$ such that
\begin{equation}\label{ratz}
g_{j}=\sum_{|i|\leq 4}r_{ji}(x)g_i  \mbox{ }\mbox{ }\mbox{ (} |j|\geq 5\mbox{)}.\end{equation} 
%If $l=4$ we get  6 equations for $i=(5,0,),(4,1),\ldots ,(0,5)$ and 6 unknowns $a_j, b_j$ for $j=(2,0),(1,1),(0,2)$.  (\ref{calcul}) 
%An important remark about this construction is the following: s
Thus  (\ref{calcul}) --  (\ref{ratz}) provide $c_{ji},r_{ji}$.  Since  $\det A  x \! \not =\! 0$ and 
$x_{40}\! =\!  \sum_{|i|=4}x_{i}t^i |_{t=e_1 }\! <\! 0$, the denominators of the rational functions $r_{ji}$ 
do not vanish at $x=\lambda^*$. $\Box$
\vspace{3 mm}

It would be interesting to generalize Proposition \ref{p1}  to arbitrary $n$ and $k$, for a  class of simple domains $T$ including $\mathbb{R}^n$, $[0,\infty )^n$  and  get rid of assumptions like $g\in G_0  , G$, for  Lagrangians $L_\epsilon$ with $\epsilon >0$. 
Also, numerical tests  of  systems like (\ref{sistem})  could be tried. 
% We hope to obtain more applications of the present results  in future work. 
%there exists a range of data $g$ for which the present approach 
%-----------------------
% defined, roughly speaking, as follows.  

Institute of Mathematics of the Czech Academy

Zitna 25, 115 67 Prague 1

Czech Republic
\vspace{1 mm}

{\it ambrozie@math.cas.cz}
\vspace{5 mm}

and: Institute of Mathematics "Simion Stoilow" - Romanian Academy, 

\hspace{8 mm} PO Box 1-764, 014700 Bucharest, Romania

\end{document}